\documentclass[preprint,3p,times,onecolumn,sort&compress]{elsarticle}

\pdfoutput=1

\usepackage{hyperref}
\usepackage{amsthm}
\usepackage{mathtools}
\usepackage{makecell}
\usepackage{booktabs}
\usepackage{multirow}
\usepackage[table,xcdraw,dvipsnames]{xcolor}
\usepackage{graphicx}
\usepackage{subcaption}
\usepackage{float}
\usepackage{overpic} 
\usepackage{rotating}
\usepackage{tikz}
\usepackage{tikz-3dplot}
\usepackage{tikzscale}
\usepackage{pgfplots}
\usetikzlibrary{calc}  
\usetikzlibrary{shapes,arrows,chains}
\usetikzlibrary{quotes,arrows.meta,angles}
\usetikzlibrary{decorations.pathreplacing}
\usetikzlibrary{positioning}
\usetikzlibrary{shapes.misc}
\usetikzlibrary{patterns}
\usetikzlibrary{3d}

\DeclareSymbolFont{rmlargesymbols}{OMX}{cmex}{m}{n}
\DeclareMathSymbol{\rmsum}{\mathop}{rmlargesymbols}{80}
\DeclareMathSymbol{\rmintop}{\mathop}{rmlargesymbols}{82}
\renewcommand{\sum}{\rmsum}
\renewcommand{\int}{\rmintop\nolimits}

\DeclareSymbolFont{greeksymbolsptm}{OML}{ptm}{m}{n}
\DeclareMathSymbol{\sigma}{\mathalpha}{greeksymbolsptm}{27}

\DeclareSymbolFont{greeksymbols}{OML}{mdbch}{m}{n}
\DeclareMathSymbol{\Omega}{\mathalpha}{greeksymbols}{10} 
\DeclareSymbolFont{greeksymbolsb}{OML}{mdbch}{b}{n}
\DeclareMathSymbol{\OM}{\mathalpha}{greeksymbolsb}{10} 		

\DeclareMathAlphabet{\altmathcal}{OMS}{cmsy}{m}{n} 	  		
\newcommand{\Sr}{\altmathcal{S}}
\newcommand{\Vr}{\altmathcal{V}}
\newcommand{\Or}{\altmathcal{O}}
\newcommand{\Qr}{\altmathcal{Q}}
\newcommand{\Hr}{\altmathcal{H}}

\renewcommand{\.}{\!\:}                             	
\newcommand{\pra}[1]{\left(#1\right)}			    	
\newcommand{\bra}[1]{\left[#1\right]}			    	
\newcommand{\ack}[1]{\left\{#1\right\}} 				
\newcommand{\prd}[2]{\left\langle#1,#2\right\rangle}	
\newcommand{\abs}[1]{\left\lvert#1\right\rvert}     	

\newcommand{\R}{\mathbb{R}}								
\newcommand{\C}{\mathbb{C}}                         	
\newcommand{\ZZ}{\mathbb{Z}}                        	
\newcommand{\ii}{\textrm{i\.}}					    	
\renewcommand{\Re}{\text{Re}\.}
\renewcommand{\Im}{\text{Im}\.}
\newcommand{\e}[1]{\text{exp}\.\big(#1\big)} 								
\newcommand{\grad}{\text{grad}}
\newcommand{\curl}{\text{curl}}
\newcommand{\curls}{\curl;\.}
\newcommand{\curlbt}{\curl^{\,\bt};\.}
\newcommand{\LOM}{(L^2(\OM))^3}
\newcommand{\LOMxz}{(L^2(\OM_{x,z}))}

\newcommand{\B}[2]{B_{#1}^{\.#2}}						
\newcommand{\x}{\xi} 									
\newcommand{\z}{\zeta}									
\newcommand{\X}{\Xi} 									
\newcommand{\Z}{\text{Z}}								
\newcommand{\nx}{n_x} 									
\newcommand{\nz}{n_z} 									
\newcommand{\px}{p_x} 									
\newcommand{\pz}{p_z} 									

\newcommand{\bF}{\textbf{F}}
\newcommand{\bE}{\textbf{E}}
\newcommand{\bM}{\textbf{M}}
\newcommand{\M}{M}
\newcommand{\bH}{\textbf{H}}
\newcommand{\ibH}{\textit{\textbf{H}}}
\renewcommand{\H}{\text{H}}
\newcommand{\Hbt}{\H^{\.\bt}}
\newcommand{\bW}{\textbf{W}}
\newcommand{\W}{\text{W}}
\newcommand{\Wbt}{\W^{\.\bt}}
\newcommand{\bn}{\mathbf{n}}

\newcommand{\bt}{\beta}
\newcommand{\om}{\omega}
\newcommand{\Om}{\Omega}
\newcommand{\sg}{\sigma}
\newcommand{\eps}{\varepsilon}

\newcommand{\az}[1]{{\color{black}#1}}
\newcommand{\fig}[1]{\mbox{\az{Fig.} #1}}

\newcommand{\eq}[1]{\mbox{Eq. #1}}
\newcommand{\eqs}[2]{\mbox{Eqs. #1 and #2}}

\newcommand{\tab}[1]{\mbox{\az{Table} #1}}

\newcommand{\Sec}[1]{\mbox{Section #1}}

\newdefinition{rem}{Remark}

\definecolor{drot}{rgb}{0.7,0,0.1}

\definecolor{C1}{rgb}{   0,    0.45,    0.74}
\definecolor{C2}{rgb}{0.85,    0.33,    0.10}
\definecolor{C4}{rgb}{   0,    0.50,    0   }
\definecolor{C3}{rgb}{0.49,    0.18,    0.56}
\definecolor{C6}{rgb}{0.47,    0.67,    0.19}
\definecolor{C5}{rgb}{0.30,    0.75,    0.93}
\definecolor{C7}{rgb}{0.64,    0.08,    0.18}
\definecolor{C8}{rgb}{1.00,       0,    1.00}
\definecolor{C9}{rgb}{0.93,    0.69,    0.13}
\definecolor{C10}{rgb}{  0,       0,    1.00}
\definecolor{C11}{rgb}{  0,       0,       0}

\pgfplotscreateplotcyclelist{mycolorlist}{%
C1, mark=none\\%
C2, mark=none\\%
C3, mark=none\\%
C4, mark=none\\%
C5, mark=none\\%
C6, mark=none\\%
C7, mark=none\\%
C8, mark=none\\%
C9, mark=none\\%
C10, mark=none\\%
C11, mark=none\\%
}

\begin{document}
\baselineskip14pt
\sloppy

\begin{frontmatter}

\title{Database Generation for Deep Learning Inversion of 2.5D Borehole Electromagnetic Measurements using Refined Isogeometric Analysis}

\author[a1]{Ali Hashemian\corref{cor1}}
\cortext[cor1]{Corresponding author}
\ead{ahashemian@bcamath.org}

\author[a4,a5]{Daniel Garcia}
\author[a2,a1]{Jon Ander Rivera}
\author[a2,a1,a3]{David Pardo}

\address[a1]{BCAM -- Basque Center for Applied Mathematics, Bilbao, Basque Country, Spain}
\address[a4]{CIMNE -- International Center for Numerical Methods in Engineering, Barcelona, Catalu\~na, Spain}
\address[a5]{IDAEA -- Institute of Environmental Assessment and Water Research, Barcelona, Catalu\~na, Spain}
\address[a2]{University of the Basque Country UPV/EHU, Leioa, Basque Country, Spain}
\address[a3]{Ikerbasque -- Basque Foundation for Sciences, Bilbao, Basque Country, Spain}

\begin{abstract}
	
Borehole resistivity measurements are routinely inverted in real-time during geosteering operations. The inversion process can be efficiently performed with the help of advanced artificial intelligence algorithms such as deep learning. These methods require a large dataset that relates multiple earth models with the corresponding borehole resistivity measurements. In here, we propose to use an advanced numerical method ---\.refined isogeometric analysis~(rIGA)\.--- to perform rapid and accurate 2.5D simulations and generate databases when considering arbitrary 2D earth models. Numerical results show that we can generate a meaningful synthetic database composed of 100,000 earth models with the corresponding measurements in 56 hours using a workstation equipped with two CPUs. 

\end{abstract}


\end{frontmatter}


\section{Introduction} 
\label{sec.Introduction}

Geosteering plays a crucial role in oil and gas engineering. 
To perform geosteering operations, companies often employ borehole instruments that record electromagnetic (EM) data in real-time while drilling~\cite{LIU2017187}.
Since the electrical resistivity is highly sensitive to salinity,
EM measurements are used to distinguish between hydrocarbon and water-saturated rocks~\cite{Liu2017}.

Inversion techniques estimate layer-by-layer EM properties from the measurements,
allowing for the adjustment of the logging trajectory during geosteering operations.
Thus, they enable to select an optimal well trajectory toward the target hydrocarbon-saturated rocks.
There exist a plethora of inversion methods in the literature, including
gradient based methods~\cite{Vogel2002,Tarantola2005}, 
statistics based methods~\cite{Watzenig2007,Kaipio2007,Shen2020}, 
and artificial intelligence based \mbox{methods~\cite{Kim2018,Yang2019,Li2020}}.
In geosteering EM measurements, deep learning (DL) methods with advanced encoder-decoder neural networks have recently demonstrated to be suitable to solve inverse problems~\cite{Shahriari2020,Shahriari20202}.

DL methods are fast, but require a massive training dataset.
To decrease the {\em online} computational time during field operations, 
we often produce such a large dataset {\em a~priori (offline)}
using tens of thousands of simulations of borehole resistivity measurements.
To generate the database for DL inversion, 
we employ simulation methods to solve Maxwell's equations with different conductivity distributions (earth models).
Since 3D simulations are expensive and possibly unaffordable when computing such large databases,
it is common to reduce the earth model dimensionality to two or one spatial dimensions
using a Fourier or a Hankel transform. 
These transformations lead to the so-called 2.5D~\cite{Abubakar2006,Pardo2008,Shen2008,Nam2013,Gernez2020}
and 1.5D~\cite{Pardo2015,Bakr2017,Shahriari20203} formulations, respectively.
1.5D simulations are inaccurate when dealing with geological faults.
In this work, we focus on the efficient generation of a database for deep learning inversion using 2.5D simulations.

Galerkin methods are effective for simulating well-logging problems (see, e.g.,~\cite{Pardo2006,Calo2011,Ma2012,Wang2013,Rodrguez2018,Chaumont2018}).
Isogeometric analysis (IGA), introduced by Hughes et al.~\cite{Hughes2005}, is a widely used Galerkin method for solving partial differential equations.
IGA has been successfully employed in various electromagnetic~\cite{Buffa2010,Nguyen2012,Buffa2014,Simpson2018,Simona2020}
and geotechnical~\cite{Shahrokhabadi2019,Hageman2019} applications.
IGA uses spline basis functions introduced in computer-aided design (CAD)
as shape functions of finite element analysis (FEA).
These basis functions exhibit high continuity (up to ${C^{p-1}}$, being $p$ the polynomial order of spline bases) across the element interfaces.

When comparing IGA and FEA, the former provides smoother solutions for wave propagation problems with a lower number of unknowns~\cite{Hughes2005,Cottrell2009}.
However, in contrast to the minimal interconnection of elements in FEA, high-continuity IGA discretizations strengthen the interconnection between elements,
leading to an increase of the cost of matrix LU factorization per degree of freedom when using sparse direct solvers~\cite{Collier2012}.
In order to avoid this degradation
and also benefit from the recursive partitioning capability of multifrontal direct solvers, 
Garcia et al.~\cite{Garcia2017} developed a new method called refined isogeometric analysis (rIGA).
This discretization technique conserves desirable properties of high-continuity IGA discretizations,
while it partitions the computational domain into blocks of macroelements weakly interconnected by low-continuity separators.
As a result, the computational cost required for performing LU factorization decreases.
The applicability of the rIGA framework to general EM problems was studied in~\cite{Garcia2019}.
Compared to high-continuity IGA,
rIGA produces solutions of EM problems up to $\Or(p^2)$ faster on large domains
and close to $\Or(p)$ faster on small domains.
rIGA also improves the approximation errors with respect to IGA since
the continuity reduction of basis functions enriches the Galerkin space.

Herein, 
we propose the use of rIGA discretizations to generate databases for DL inversion of
2.5D geosteering EM measurements.
We consider {\em a~priori} grids following the idea of optimal grid generation for 2.5D EM measurements presented by Rodr\'iguez-Rozas et al.~\cite{Rodrguez2018}
and the methods described in~\cite{Rodrguez2016} for Fourier mode selections.
Compared to the FEA approach described in~\cite{Rodrguez2018} that assigns increasing polynomial orders for the elements near the well, 
we consider a (smooth) high-continuity IGA discretization with a fixed polynomial order everywhere and reduce the computational cost by continuity reduction of certain basis functions in the sense of rIGA framework.
To assess the accuracy and computational efficiency of the rIGA approach in borehole resistivity simulations,
we consider several model problems where high-angle wells cross spatially heterogeneous media exhibiting multiple geological faults.
Then, we investigate the performance of the proposed approach when 
generating a synthetic database composed of 100,000 earth models for DL inversion.

The remainder of this article is organized as follows.
In \Sec{\ref{sec.2.5D}}, we review the governing equations of the 3D EM wave propagation problem and derive the 2.5D variational formulation.
\Sec{\ref{sec:ProblemDescription}} introduces the considered borehole resistivity problem.
In \Sec{\ref{sec.IGArIGA}}, we describe both high-continuity and refined isogeometric discretizations,
followed by the implementation details in \Sec{\ref{sec.Implementation}}.
We analyze the accuracy and computational efficiency of the rIGA approach when applied to borehole resistivity problems in \Sec{\ref{sec.Results}}.
In this section, we also generate a database for DL inversion of geosteering measurements.
Finally, \Sec{\ref{sec.Conclusions}} draws the main conclusions
and possible future research lines stemming from this work.


\section{2.5D variational formulation of EM measurements} 
\label{sec.2.5D}

\subsection{3D wave propagation problem} 
\label{sub:3D}

The two time-harmonic curl Maxwell's equations describing the 3D wave propagation in an isotropic medium are

\begin{align}
	\nabla\times\bE + \ii\om\mu\bH &= -\ii\om\mu\bM \,, \\
	\nabla\times\bH 			         &= (\sg+\ii\om\eps)\bE \,,
\end{align}	

\noindent
where $\bE$ is the electric field, $\bH$ is the magnetic field, $\ii$ is the imaginary unit, $\sg$ is the electric conductivity, $\eps$ is the electric permittivity, $\mu$ is the magnetic permeability, ${\om=2\pi f}$ is the angular frequency, being $f$ the source frequency, and $\bM$ is the time-harmonic magnetic source located at ${(x_0,y_0,z_0)}$ and given by
\begin{align}
	\bM = \delta(x-x_0)\.\delta(y-y_0)\.\delta(z-z_0)\big[\M_x\,,\M_y\,,\M_z\big]^T \quad {\rm in}~\R^3.
\end{align}

From Maxwell's equations, we obtain the following reduced wave formulation in terms of magnetic field $\bH$:

\begin{equation}
	\begin{cases}
	\begin{aligned}
		&\text{Find }\bH=\big[\H_x\,,\H_y\,,\H_z\big]^T, \text{ with } \bH:\OM\subset\R^3\rightarrow\C^3, \text{ such that:} \\
		&\qquad\begin{aligned}
      \nabla\times\pra{\dfrac{1}{\sg+\ii\om\eps}\nabla\times\bH} + \ii\om\mu\bH &= -\ii\om\mu\bM \quad&&\text{in }\OM \,,\\
		  \bE\times\bn &= \mathbf{0} &&\text{on }\partial{\OM} \,,\end{aligned}
	\end{aligned}		
	\end{cases}
	\label{eq:str}
\end{equation}

\noindent
where $\OM$ is the domain of study defined as a tensor-product box by 

\begin{equation}
\begin{aligned}
	\OM &= \Om_x\times\Om_y\times\Om_z \\
	   	&= \big(-L_x/2,L_x/2\big) \times \big(-L_y/2,L_y/2\big) \times \big(-L_z/2,L_z/2\big) \,,
\end{aligned}		
\end{equation}

\noindent
being $L_x$,$L_y$, and $L_z$ positive real constants, and ${\partial{\OM}}$ the domain boundary.

To introduce the weak formulation of this problem, we first define the $\ibH(\curls\OM)$-conforming functional spaces 
\begin{align}
	\ibH(\curls\OM)      &= \ack{\bW=\big[\W_x\,,\W_y\,,\W_z\big]^T\in\LOM:\nabla\times\bW\in\LOM} \,,\\
	\ibH_{0}(\curls\OM)  &= \Big\{\bW\in\ibH(\curls\OM):\bW\times\bn=\mathbf{0}\text{ on }\partial{\OM}\Big\} \,.
\end{align}	

\noindent 
The $\ibH(\curls\OM)$ space is endowed with the inner product

\begin{equation}
\begin{aligned}
	(\bW,\bH)_{\ibH(\curls\OM)} &\coloneqq (\nabla\times\bW,\nabla\times\bH)_{\LOM} + (\bW,\bH)_{\LOM} \\ 
	&\coloneqq \int_{\OM} \prd{\nabla\times\bW}{\nabla\times\bH}\,d\OM + \int_{\OM} \prd{\bW}{\bH}\,d\OM \,,
\end{aligned}
\end{equation}

\noindent 
where $\prd{\cdot\,}{\cdot}$ denotes the Hermitian $L^2$ inner product on the complex vector space. 

We build the weak formulation by multiplying \eq{\eqref{eq:str}} with an arbitrary function ${\bW \in \ibH_{0}(\curls\OM)}$, using Green's formula, and integrating over the domain $\OM$. The weak formulation is then

\begin{equation}
	\begin{cases}
	\begin{aligned}
		&\text{Find }\bH\in\ibH_{0}(\curls\OM), \text{ such that for every } \bW\in\ibH_{0}(\curls\OM) \,, \\
		&\qquad\pra{\nabla\times\bW,\dfrac{1}{\sg+\ii\om\eps}\nabla\times\bH}_{\LOM} + \ii\om\mu\.(\bW,\bH)_{\LOM} = -\ii\om\mu\.(\bW,\bM)_{\LOM}\,.
	\end{aligned}		
	\end{cases}
	\label{eq:wke}
\end{equation}


\subsection{2.5D variational formulation} 
\label{sub.2.5D}

Herein, we focus on the case when the material properties are homogeneous along one spatial direction, e.g., $y$-axis. 
We denote the domain for this case as ${\OM\coloneqq\Om_y\times\OM_{x,z}}$.
We perform a Fourier transform along the $y$-axis to represent the 3D problem as a sequence of uncoupled 2D problems, one per Fourier mode.
In this case, we define the magnetic field $\bH$ as a series expansion using the complex exponentials: 
\begin{align}
	\bH \coloneqq \sum_{\bt=-\infty}^{+\infty} \bH_{\bt}\,\e{\ii 2\pi\bt y/L_y} \,.
  \label{eq:FourierSeries}
\end{align}

\noindent
being $\bt$ the Fourier mode number and ${\bH_{\bt}=\big[\Hbt_x\,,\Hbt_y\,,\Hbt_z \big]^T}$ with ${\bH_{\bt}:\OM_{x,z}\subset\R^2\rightarrow\C^3}$.
Fourier modes satisfy the following orthogonality relationships: 
\begin{align}
  \dfrac{1}{L_y}\int_{-L_y/2}^{L_y/2} \e{\ii 2\pi\bt_1 y/L_y}\,\e{\ii 2\pi\bt_2 y/L_y}\,dy = \delta_{\bt_1}\delta_{\bt_2} \,.
\end{align}

By employing a test function of the form

\begin{equation}
	\bW \coloneqq \dfrac{1}{L_y}\bW_{\bt}\,\e{\ii 2\pi\bt y/L_y} \,,
\end{equation}

\noindent
and using the $\ibH(\curlbt\OM_{x,z})$-conforming functional spaces
\begin{align}
  \ibH(\curlbt\OM_{x,z}) &= \ack{\bW_{\bt}=\big[\Wbt_x\,,\Wbt_y\,,\Wbt_z\big]^T\in\LOMxz^3:\Wbt_y\in H^1(\OM_{x,z})\text{ and }\nabla\times\big[\Wbt_x\,,\Wbt_z\big]^T\in\LOMxz^2} \,,\\
  \ibH_{0}(\curlbt\OM_{x,z}) &= \Big\{\bW_{\bt}\in\ibH(\curlbt\OM_{x,z}):\bW_{\bt}\times\bn=\mathbf{0}\text{ on }\partial{\OM} \Big\} \,,
\end{align}

\noindent
we build the following variational formulation from \eq{\eqref{eq:str}} by integrating over $\OM_{x,z}$:

\begin{equation}
	\begin{cases}
	\begin{aligned}
		&\text{Find } \bH=\textstyle\sum_{\bt=-\infty}^{+\infty} \bH_{\bt}\,\e{\ii 2\pi\bt y/L_y} \,, ~ \bH_{\bt}\in\ibH_{0}(\curlbt\OM_{x,z})\\
		&\text{such that for every }\bt\in\ZZ\text{ and }\bW_{\bt}\in\ibH_{0}(\curlbt\OM_{x,z}) \,,\\
		&\qquad\pra{\nabla^{\bt}\times\bW_{\bt},\dfrac{1}{\sg+\ii\om\eps}\nabla^{\bt}\times\bH_{\bt}}_{\LOMxz^3} + \ii\om\mu\pra{\bW_{\bt},\bH_{\bt}}_{\LOMxz^3} = -\ii\om\mu\pra{\bW_{\bt},\bM_{\bt}}_{\LOMxz^3} \,,
	\end{aligned}		
	\end{cases}
	\label{eq:2.5Dwke}
\end{equation}

\noindent where 
\begin{align}
	\nabla^{\bt}\times\bW_{\bt} \coloneqq \bra{\ii\bt\.\dfrac{2\pi}{L_y}\Wbt_z-\dfrac{d\Wbt_y}{dz}\,,~ \dfrac{d\Wbt_x}{dz}-\dfrac{d\Wbt_z}{dx}\,,~ \dfrac{d\Wbt_y}{dx}-\ii\bt\.\dfrac{2\pi}{L_y}\Wbt_z}^T ,
\end{align}

\noindent
and $\bM_{\bt}$ is the time-harmonic magnetic source written in terms of the Fourier transform as
\begin{align}
	\bM_{\bt} = \dfrac{1}{L_y}\.\delta(x-x_0)\.\delta(z-z_0)\big[\M_x\,,\M_y\,,\M_z\big]^T \e{\ii 2\pi\bt y_0/L_y} \,.
\end{align}

This formulation corresponds to the 2.5D variational formulation 
previously described in, e.g.,~\cite{Rodrguez2018,Chaumont2018}.


\section{Borehole resistivity measurement acquisition system} 
\label{sec:ProblemDescription}

\begin{figure}[!h]
\centering
\includegraphics{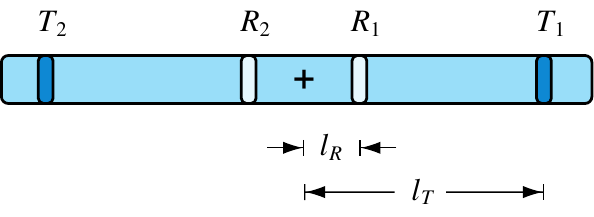}
	\caption{A schematic LWD instrument with two transmitters and two receivers symmetrically located around the tool center.}
\label{fg:LoggingTool}
\end{figure}

We consider a logging-while-drilling (LWD) instrument
equipped with transmitters ($T_i$) and receivers ($R_j$).
This tool is sensitive to resistivities within the range ${0.2\sim500~\Om\cdot\rm m}$ (phase resistivity) and ${0.2~\sim300~\Om\cdot\rm m}$ (amplitude resistivity) under an operating frequency between ${\rm 0.1~and~2~MHz}$~\cite{LIU2017187}.
For the sake of simplicity, herein, we restrict to two transmitters and two receivers symmetrically located around the tool center (see \fig{\ref{fg:LoggingTool}})
at an operating frequency of 2 MHz.

Triaxial logging instruments generate measurements for all possible orientations of the transmitter--receiver pairs.
We follow the notation presented in~\cite{Davydycheva2011,Rodrguez2018} to denote the magnetic field.
Thus, we write ${\H_{ZZ}^{T_iR_j}\in\C}$ as the coaxial magnetic field in the borehole system of coordinates induced by source $T_i$ and measured at receiver $R_j$ (${i,j=1,2}$).
We use the magnetic fields measured at $R_1$ and $R_2$ to compute the {\em attenuation ratio} and {\em phase difference}.
We symmetrize the signal originating from $T_1$ and $T_2$ to obtain the {\em quantity of interest} $\text{Q}_{ZZ}$ at each logging position as
\begin{align}
  \text{Q}_{ZZ} \coloneqq \dfrac{1}{2}\pra{
  \log\dfrac{\H_{ZZ}^{T_1R_1}}{\H_{ZZ}^{T_1R_2}} + \log\dfrac{\H_{ZZ}^{T_2R_2}}{\H_{ZZ}^{T_2R_1}}} \,.
\end{align}

\noindent
Then, we compute the attenuation ratio $(A)$ and phase difference $(P)$, respectively, as the real and imaginary parts of~$\text{Q}_{ZZ}$:
\begin{align}
  A &\coloneqq \Re(\text{Q}_{ZZ}) \label{eq:Real} \,,\\
  P &\coloneqq \Im(\text{Q}_{ZZ}) \label{eq:Image} \,.
\end{align}

\noindent
We can then obtain the {\em apparent resistivities} based on the attenuation and phase (\,${\rho_A}$ and ${\rho_P}$, respectively) using a {\em look-up} algorithm presented in~\cite{Anderson2001}.


\section{Refined isogeometric analysis} 
\label{sec.IGArIGA}

In this work, we consider a {\em multi-field} EM problem and 
discretize the 2.5D variational formulation of \eq{\eqref{eq:2.5Dwke}} using a B-spline generalization of a curl-conforming space, introduced by Buffa et al.~\cite{Buffa2010}.
We first review some basic concepts of high-continuity IGA discretizations.


\subsection{High-continuity IGA discretization} 
\label{sub:IGA}

Given the parametric domain ${\big\{\x,\z\in\hat{\OM}_{x,z}:(0,1)^2\subset\R^{2}\big\}}$,
we introduce the spline space ${\Sr^{\px,\pz}_{k_x,k_z}}$ as

\begin{equation}
	\Sr^{\px,\pz}_{k_x,k_z} \coloneqq \text{span} \ack{\B{i,j}{\px,\pz}}_{i=0,j=0}^{\nx-1,\.\nz-1} \,,
\end{equation}

\noindent
where $n$, $p$, and $k$ with their indices are the number of degrees of freedom, polynomial degree, and continuity of basis functions in $x$ and $z$ directions, respectively.
The bivariate basis functions are

\begin{equation}
	\B{i,j}{\px,\pz} \coloneqq \B{i}{\px}(\x) \otimes \B{j}{\pz}(\z) \,, \qquad i=0,1,...,\nx-1 \,, \quad j=0,1,...,\nz-1 \,,
\end{equation}

\noindent
where the univariate bases are expressed by the Cox--De Boor recursion formula~\cite{TheNURBSBook} as 

\begin{equation}
\begin{aligned}
	\B{i}{0}(\x)  &= \left\{\begin{array}{ll} 
		1 & \x_i\leq \x < \x_{i+1}\,, \\ 
		0 & \mbox{otherwise}\,, \end{array} \right. \\
	\B{i}{p}(\x) &= \dfrac{\x-\x_i}{\x_{i+p}-\x_i}\B{i}{{p-1}}(\x)+
	\dfrac{\x_{i+p+1}-\x}{\x_{i+p+1}-\x_{i+1}}\B{{i+1}}{{p-1}}(\x)\,,
\end{aligned}
\end{equation}

\noindent
and spanned over the respective knot sequences in $x$ and $z$ directions, given by
\begin{align}
	\X & =[\underbrace{0,0,...,0}_{\px+1},\x_{\px+1},\x_{\px+2}...,\x_{\nx-1},\underbrace{1,1,...,1}_{\px+1}]\,, \\
	\Z & =[\underbrace{0,0,...,0}_{\pz+1},\z_{\pz+1},\z_{\pz+2}...,\z_{\nz-1},\underbrace{1,1,...,1}_{\pz+1}]\,.
\end{align}

\noindent
We assume single multiplicities for all knots, providing maximum continuity $k=p-1$ for the IGA discretization.

\begin{figure}[!h]
	\centering
    \includegraphics{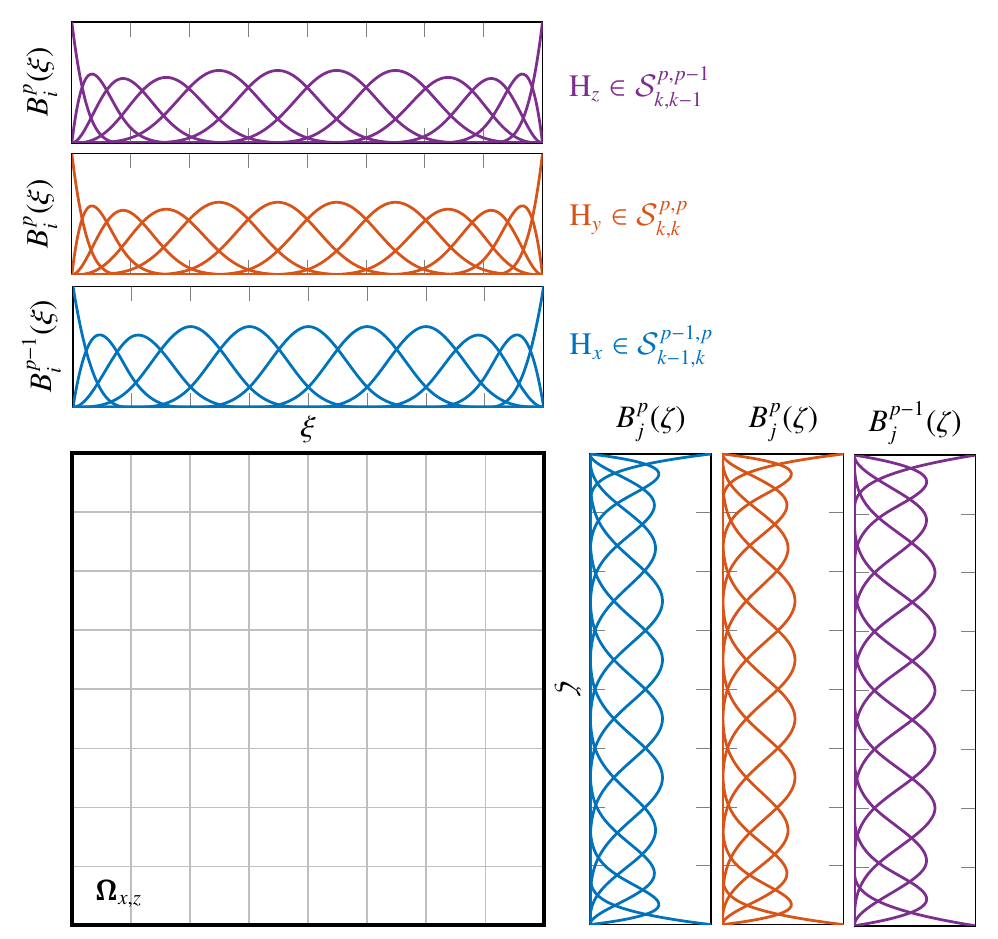}
	\caption{Example of the ${\ibH(\curl)\times H^1}$ space for a 2.5D formulation discretized by $C^{p-1}$ IGA with uniform ${8\times8}$ elements in ${\OM_{x,z}}$, polynomial degree ${p=4}$, and continuity ${k=3}$. The univariate basis functions of $\H_x$, $\H_y$, and $\H_z$ are shown in blue, red, and purple, respectively. Thin gray lines in the mesh skeleton denote the high-continuity element interfaces.}
	\label{fg:IGA_curlh1}
\end{figure}

\fig{\ref{fg:IGA_curlh1}} illustrates the ${\ibH(\curl)\times H^1}$ IGA discrete space in $\OM_{x,z}$
along with the univariate basis functions of the respective vector fields.
For brevity, herein and in the following, we exclude the superscript $\bt$ in referring to the components of the magnetic field.

We define the spaces in the parametric domain and introduce the appropriate transformations to obtain the discretization on the physical domain. 
We start with the set of discrete spaces in the parametric domain, given by
\begin{align}
	\hat{\Vr}_h^{\,\curl}(\hat{\OM}_{x,z}) &\coloneqq \Sr^{p-1,p}_{k-1,k} \times \Sr^{p,p-1}_{k,k-1} \,, \\
	\hat{\Qr}_h^{\,\grad}(\hat{\OM}_{x,z}) &\coloneqq \Sr^{p,p}_{k,k} \,.
\end{align}

\noindent
By defining ${\bF:\hat{\OM}_{x,z}\rightarrow\OM_{x,z}}$ as the geometric mapping from the parametric domain onto the physical domain, and $D\bF$ as its Jacobian, 
we introduce the set of discrete spaces in the physical domain:
\begin{align}
 	\Vr_{h}^{\,\curl}(\OM_{x,z}) &\coloneqq \ack{\bH_{x,z}=\big[\H_x\,,\H_z\big]^T\in\ibH(\curlbt\OM_{x,z})\cap\ibH^1(\OM_{x,z}):\iota^{\.\curl}(\bH_{x,z})=\hat{\bH}_{x,z}\in\hat{\Vr}_{h}^{\,\curl}(\hat{\OM}_{x,z})} \,, \\
 	\Qr_{h}^{\,\grad}(\OM_{x,z}) &\coloneqq \ack{\H_y\in H^1(\OM_{x,z}):\iota^{\.\grad}(\H_y)=\hat{\H}_y\in\hat{\Qr}_h^{\,\grad}(\hat{\OM}_{x,z})} \,, 
 	\label{eq:map_curl}
\end{align}

\noindent
where
we use the following curl- and grad-preserving pullback mappings~\cite{Garcia2019,Buffa2010}:
\begin{align}
	\iota^{\.\curl}(\bH_{x,z}) 	&\coloneqq (D\bF)^{T}(\bH_{x,z} \circ \bF) \,, \\
	\iota^{\.\grad}(\H_y)	&\coloneqq \H_y \circ \bF\,.
\end{align}

Thus, by defining the discrete space
\begin{align}
\Hr_{h,0}(\OM_{x,z}) \coloneqq \ack{\bH_{\bt,h}\in\Vr_{h}^{\,\curl}(\OM_{x,z}) \times \Qr_{h}^{\,\grad}(\OM_{x,z}) : \bH_{\bt,h}\times\bn=\mathbf{0}\text{ on }\partial{\OM} } \,, 
\end{align}

\noindent
we write the discrete form of \eq{\eqref{eq:2.5Dwke}} as follows (subscript $h$ refers to discrete solution):
\begin{align}
\begin{cases}
\begin{aligned}
&\text{Find } \bH_h=\textstyle\sum_{\bt=-\infty}^{+\infty} \bH_{\bt,h}\,\e{\ii 2\pi\bt y/L_y} \,, ~\bH_{\bt,h}\in \Hr_{h,0}(\OM_{x,z}) \\[3pt]
&\text{such that for every }\bt\in\ZZ {\rm~and~} \bW_{\bt,h}\in \Hr_{h,0}(\OM_{x,z}) \,,\\[2pt]
&\qquad\pra{\nabla^{\bt}\times\bW_{\bt,h}\,,\dfrac{1}{\sg+\ii\om\eps}\nabla^{\bt}\times\bH_{\bt,h}}_{\LOMxz^3} + \ii\om\mu\pra{\bW_{\bt,h},\bH_{\bt,h}}_{\LOMxz^3} = -\ii\om\mu\pra{\bW_{\bt,h},\bM_{\bt}}_{\LOMxz^3} \,.
\end{aligned}
\end{cases}
\label{eq:2.5discrete}
\end{align}


\subsection{rIGA discretization} 
\label{sub:rIGA}

The refined isogeometric analysis (rIGA) is a discretization technique
that optimizes the performance of direct solvers.
In particular, rIGA preserves the optimal convergence order of the direct solvers with respect to a fixed number of elements in the domain. 
Garcia et al.~\cite{Garcia2017} first presented this strategy for $H^1$ spaces
and then extended it to $\ibH(\curl)$, $\ibH(\text{div})$, and $L^2$ spaces (see~\cite{Garcia2019}).
Starting from the high-continuity $C^{p-1}$ IGA discretization, rIGA reduces the continuity of certain basis functions 
by increasing the multiplicity of the respective existing knots.
Hence, the computational domain is subdivided into high-continuity macroelements interconnected by low-continuity hyperplanes.
These hyperplanes coincide with the locations of the {\em separators} at different partitioning levels of the multifrontal direct solvers. 
Thus, rIGA reduces the computational cost of matrix factorization when solving PDE systems in comparison to IGA and FEA.

\begin{figure}[!pb]
	\centering
    \includegraphics{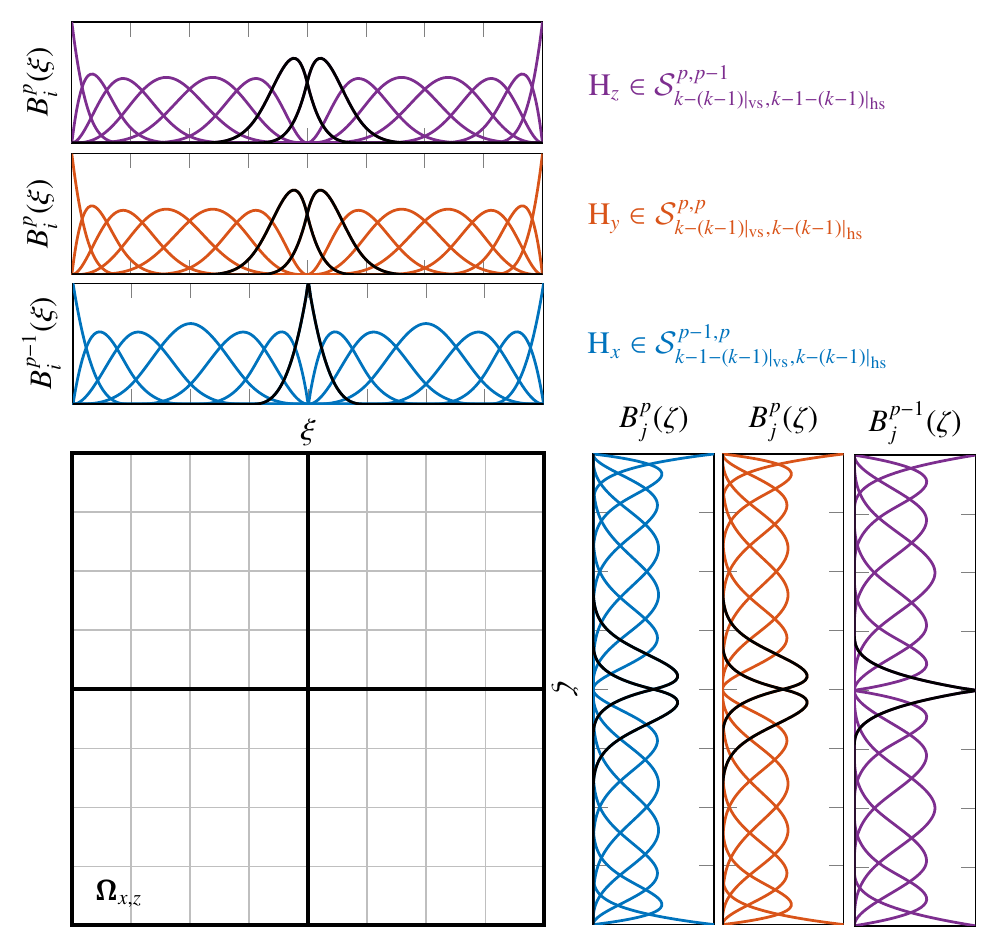}
	\caption{${\ibH(\curl)\times H^1}$ rIGA space of the ${8\times 8}$ domain of \fig{\ref{fg:IGA_curlh1}} with ${p=4}$ and ${k=3}$,
	after one level of symmetric partitioning by the rIGA discretization that results in ${4\times4}$ macroelements.
	rIGA reduces the continuity of basis functions by $k-1$ degrees across the macroelement separators
	(the low-continuity bases are shown in black).
	Thin gray lines in the mesh skeleton denote the high-continuity element interfaces, 
	while thick black lines illustrate the macroelement boundaries.
	We refer to the vertical and horizontal separators as ``vs" and ``hs", respectively.}
	\label{fg:rIGA_curlh1}
\end{figure}

For multi-field problems discretized using ${\ibH(\curl)\times H^1}$ spaces, 
we preserve the commutativity of the \mbox{\em de Rham} diagram~\cite{Demkowicz2000} by reducing the continuity in ${k-1}$ degrees. 
To achieve this, we use both $C^{\.0}$ and $C^{1}$ hyperplanes and reduce the continuity across the interface between the subdomains (i.e., macroelements).
\fig{\ref{fg:rIGA_curlh1}} depicts the rIGA discretization of the ${\ibH(\curl)\times H^1}$ space of \fig{\ref{fg:IGA_curlh1}} after one level of symmetric partitioning, 
which results in macroelements containing ${4\times4}$ elements.

Previous works show rIGA discretizations provide large improvements in the solution time and memory requirements.
In particular, 
the rIGA solution is obtained up to $\Or(p^2)$ faster in large domains 
---\.and $\Or(p)$ faster in small domains\.---
than the IGA solution.
In comparison to traditional FEA with the same number of elements, rIGA provides even larger improvements.
rIGA also reduces the memory requirements since the rIGA LU factors have fewer nonzero entries than the IGA LU factors.
Finally, rIGA improves the approximation error with respect to IGA since
the continuity reduction of basis functions enriches the Galerkin space (see~\cite{Garcia2017,Garcia2019}).


\section{Implementation details} 
\label{sec.Implementation}

We implement discrete ${\ibH(\curl)\times H^1}$ spaces using PetIGA-MF~\cite{petigamf}, a multi-field extension of PetIGA~\cite{petiga}, a high-performance isogeometric analysis implementation based on PETSc (portable extensible toolkit for scientific computation)~\cite{petsc}.
PetIGA-MF allows the use of different spaces for each field of interest
and employs data management libraries to condense the data of multiple fields in a single object, thus simplifying the discretization construction.
This framework also allows us to investigate both IGA and rIGA discretizations in our 2.5D problem
with different numbers of elements, different polynomial degrees of the B-spline spaces, and different partitioning levels of the mesh.

We use Intel MKL with PARDISO~\cite{pardiso1,pardiso2} as our sparse direct solver package to construct LU factors for solving the linear systems of equations.
PARDISO employs supernode techniques to perform the matrix factorization (see, e.g.,~\cite{Schenk2000,Schenk2004}).
It provides parallel factorization using OpenMP directives~\cite{Dagum1998}
and uses the automatic matrix reordering provided by METIS~\cite{Karypis1998}.
We executed all tests on a workstation
equipped with two Intel Xeon Gold 6230 CPUs at 2.10 GHz
with 40 threads per CPU.

\begin{figure}[!h]
\centering
\includegraphics{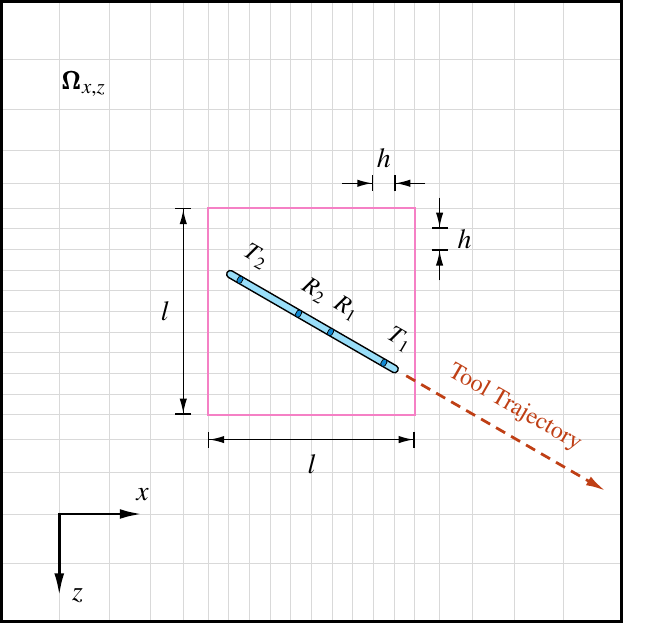}
\caption{A drawing of the computational domain $\OM_{x,z}$ and the tool trajectory.
The central subgrid bounded by a magenta box is composed of a set of fine elements located in the proximity of the logging instrument.
The remaining elements grow smoothly in size until reaching the boundary.}
\label{fg:mesh}
\end{figure}

We employ a tensor-product mesh with variable element sizes (see \fig{\ref{fg:mesh}}).
At each logging position, the computational mesh has a fine subgrid in the central part of the domain with element size equal to ${h \times h}$. 
This subgrid is surrounded by another tensor-product grid whose element sizes grow slowly until reaching the boundary.\\

\begin{rem}
For each logging position, we perform a single {\em symbolic} factorization common to all Fourier modes, followed by a numerical factorization per Fourier mode.
Once we solve the system of equations for the first transmitter, we update the right-hand side of \eq{\eqref{eq:2.5discrete}} to solve for the magnetic field induced by the second transmitter and use backward substitution.
Hence, we perform only one LU factorization per Fourier mode per logging position for both transmitters.
\end{rem}

\begin{rem}
The convergence of the Fourier series leads to a fast decay of the real and imaginary parts of $\H_{ZZ}$ for higher Fourier modes (see~\cite{Rodrguez2016}).
Thus, we truncate the series of \eq{\eqref{eq:FourierSeries}} when the magnetic field at the receivers is sufficiently small, 
such that ${\bt\in[-N_f,N_f]}$, being $2N_f+1$ the total number of Fourier modes.
In here, because of the symmetry of the media along the $y$ direction, we only consider ${\bt\in[0,N_f]}$.
\end{rem}


\section{Numerical Results}
\label{sec.Results}

In this section,
we first assess the accuracy of the rIGA approach in 
a homogeneous medium.
We also investigate the computational efficiency of the rIGA framework in comparison with IGA and FEA approaches.
Then, 
we consider two model problems consisting of high-angle wells crossing
spatially heterogeneous media 
with multiple geological faults.
Finally,
we produce our synthetic training dataset for DL inversion.
In our simulations,
we consider one operational mode of a commercial logging tool~\cite{Zhou2016}
with ${l_{\.R}=10.16~\rm cm}$ and ${l_{\.T}=56.8325~\rm cm}$ (see \fig{\ref{fg:LoggingTool}}).
We select the free space electric permittivity and magnetic permeability, i.e.,
${\eps=8.85\times10^{-12}~\rm F\cdot m^{-1}}$ and
${\mu=4\pi\times10^{-7}~\rm N\cdot A^{-2}}$, respectively.
We also consider a source frequency ${f=2~\rm MHz}$.


\subsection{Homogeneous medium}
\label{sub:HomogeneousMedia}

We assume the logging instrument is placed in a homogeneous medium with 
electrical conductivity ${\sg=0.01~(\Om\cdot\rm m)^{-1}}$.
This high-resistivity case is numerically more challenging than low-resistivity cases since it requires a larger number of Fourier modes and numerical precision.


\subsubsection{Accuracy assessment} 
\label{ssub:accuracy}

\begin{figure}[!h]
  \centering
  \begin{subfigure}{0.49\textwidth}\centering
    \pgfplotsset{width=1.0\textwidth}
    \includegraphics{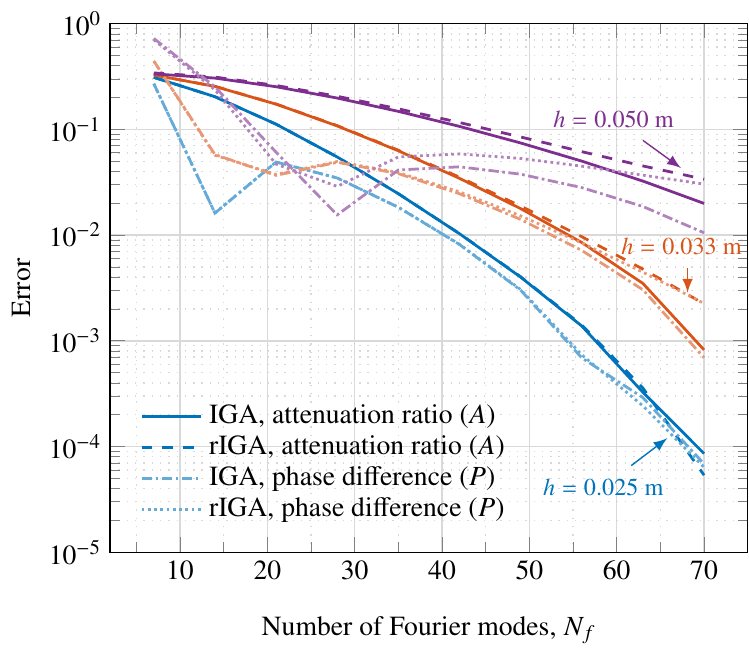}
    \caption{$p=3$}
    \label{fg:hbeta.p3}
  \end{subfigure}
  \begin{subfigure}{0.49\textwidth}\centering
    \pgfplotsset{width=1.0\textwidth}
    \includegraphics{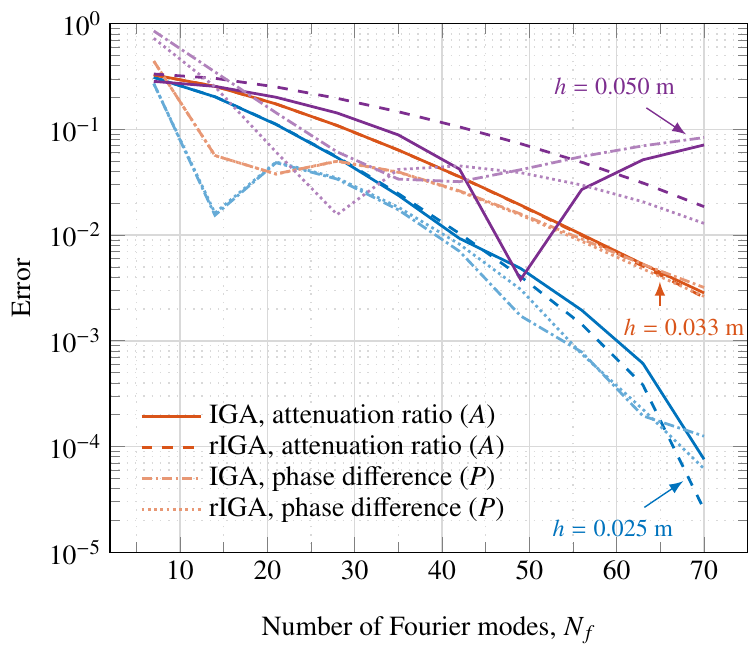}
    \caption{$p=4$}
    \label{fg:hbeta.p4}
  \end{subfigure}  
  \begin{subfigure}{0.49\textwidth}\centering
    \pgfplotsset{width=1.0\textwidth}
    \includegraphics{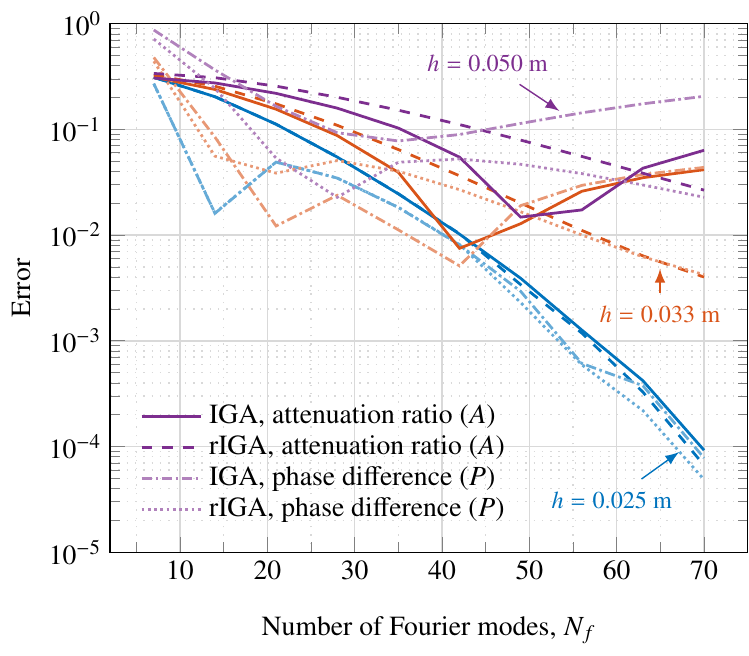}
    \caption{$p=5$}
    \label{fg:hbeta.p5}
  \end{subfigure}   
  \caption{Numerical errors when computing attenuation ratio, $A$, and phase difference, $P$, in a homogeneous medium using IGA and rIGA discretizations, obtained by for a ${64\times64}$ mesh with different element sizes $h$ and polynomial degrees $p$.}
  \label{fg:hbeta}
\end{figure}

To assess the accuracy and select certain discretization parameters,
we compare the numerical attenuation ratio, $A$, and phase difference, $P$, given by \eqs{\eqref{eq:Real}}{\eqref{eq:Image}},
with the expected (i.e., exact) values, $A_e$ and $P_e$, obtained from \.${\rho_e=1/\sg}$.
\fig{\ref{fg:hbeta}} shows the numerical errors, 
i.e., ${\abs{1-A/A_e}}$ and ${\abs{1-P/P_e}}$,
as a function of the number of Fourier modes
when computing attenuation and phase in the homogeneous medium.
Herein, we select a domain with ${64\times64}$ elements to ensure a fast numerical solution for our measurements.
We compare the high-continuity ${C^{p-1}}$ IGA with an rIGA discretization that employs ${8\times8}$ macroelements. 
This macroelement size provides the fastest results for moderate size domains (see~\cite{Garcia2017}). 
We consider three different mesh sizes
---\.${h=0.025~{\rm m},~0.33~{\rm m},~{\rm and}~0.50~{\rm m}}$\.---
and different polynomial degrees
---\.${p=3,~4,{~\rm and}~5}$.
The best results correspond to $h=0.025~{\rm m}$ (blue lines in the figure).
We also observe that rIGA discretizations (dashed and dotted lines) deliver lower errors compared to their IGA counterparts.

To investigate the decay of the solution for each Fourier mode, we compare numerical results with
the analytical 2.5D solution in the homogeneous medium presented in~\cite{Rodrguez2018}.
In particular, given $\M_z$ as the only nonzero component of the magnetic source, it is possible to analytically determine the coaxial magnetic field for each Fourier mode as
\begin{align}
  \H_{ZZ}(\,\bt\.) = -\ii\om\mu\sg\tau_{\bt_z}+\dfrac{\partial^{\.2}\tau_{\bt_z}}{\partial z^2} \,,
\end{align}

\noindent with
\begin{align}
  \tau_{\bt_z} = \dfrac{\M_z}{2\pi}\dfrac{1}{L_y}K_0({\rm C}R)\,\e{\ii2\pi\bt y_0/L_y} \,,
\end{align}

\noindent where $K_0(\cdot)$ is the modified Bessel function of the second kind of order zero, and
\begin{align}
  {\rm C} &= \pra{2\pi\bt/L_y}^2+\ii\om\mu\sg \,, \\
  R &= \sqrt{(x-x_0)^2+(z-z_0)^2} \,.
\end{align}

\begin{figure}[!h]
  \centering
  \begin{subfigure}{0.49\textwidth}\centering
    \pgfplotsset{width=1.0\textwidth}
    \includegraphics{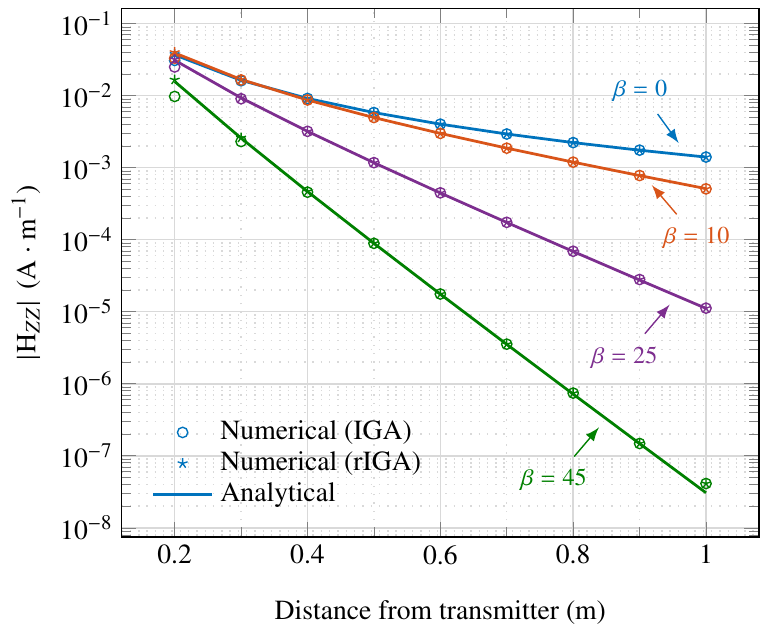}
    \caption{$p=3$}
    \label{fg:decay.p3}
  \end{subfigure}
  \begin{subfigure}{0.49\textwidth}\centering
    \pgfplotsset{width=1.0\textwidth}
    \includegraphics{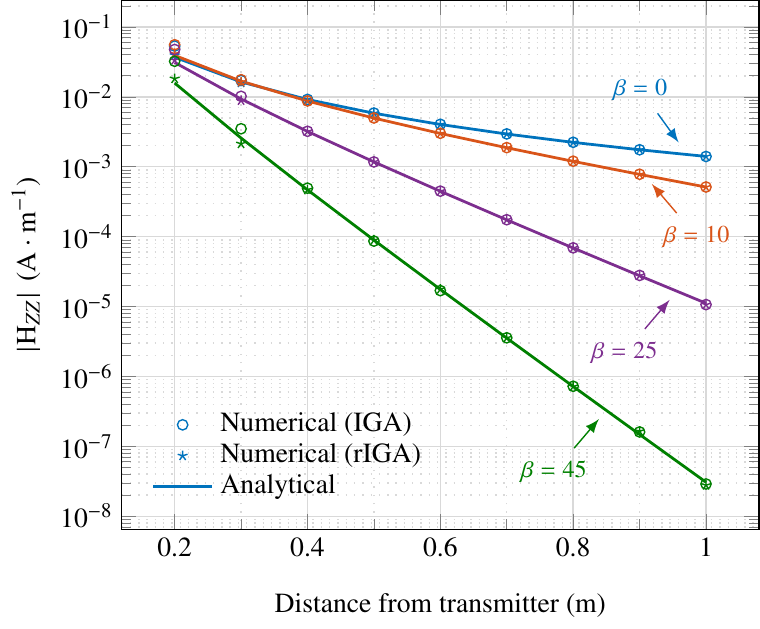}
    \caption{$p=4$}
    \label{fg:decay.p4}
  \end{subfigure} \\[6pt]
  \begin{subfigure}{0.49\textwidth}\centering
    \pgfplotsset{width=1.0\textwidth}
    \includegraphics{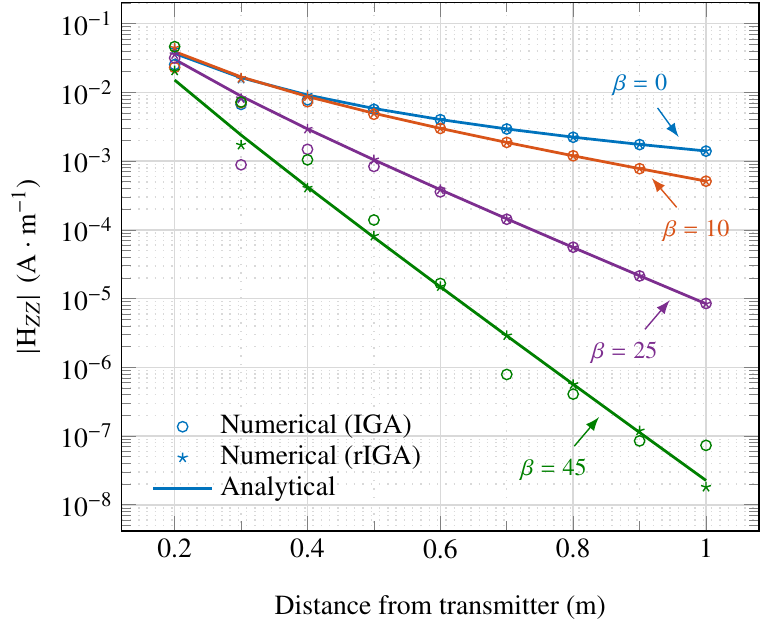}
    \caption{$p=5$}
    \label{fg:decay.p5}
  \end{subfigure}
  \caption{Comparison of the decay of the numerical and analytical coaxial magnetic fields for some Fourier modes, obtained in a grid of ${64\times64}$ elements with ${h=0.025~{\rm m}}$ and different polynomial degrees.}
  \label{fg:decay}
\end{figure}

\fig{\ref{fg:decay}} compares the decay
of the numerical coaxial magnetic field ${\H_{ZZ}(\,\bt\.)}$ 
with its analytical counterpart for some Fourier modes.
Using a domain with ${64\times64}$ elements and ${h=0.025~{\rm m}}$, we monitor the decay of the propagated waves at distances within the interval ${[0.2,1.0]~{\rm m}}$ from the transmitters
to ensure that the solutions at both receivers properly approximate the analytical ones.
Results show that rIGA discretizations deliver increased accuracy for all tested polynomial degrees.

\subsubsection{Computational efficiency} 
\label{subsub:efficiency}

Works~\cite{Garcia2017,Garcia2019} provide theoretical cost estimates of solving $H^1$ and $\ibH(\curl)$ discrete spaces, respectively. 
Herein, we add these estimates to predict the cost of discretizing the
${\ibH(\curl)\times H^1}$ space appearing in our 2.5D EM problem. 
We conclude that the cost of LU factorization of the rIGA matrix for this combined space is between ${\Or(p)}$ and ${\Or(p^2)}$ times smaller than that for IGA. Details are omitted for the sake of simplicity.

\begin{figure}[!h]
  \centering
  \begin{subfigure}{0.49\textwidth}\centering
    \pgfplotsset{height=0.85\textwidth}
    \includegraphics{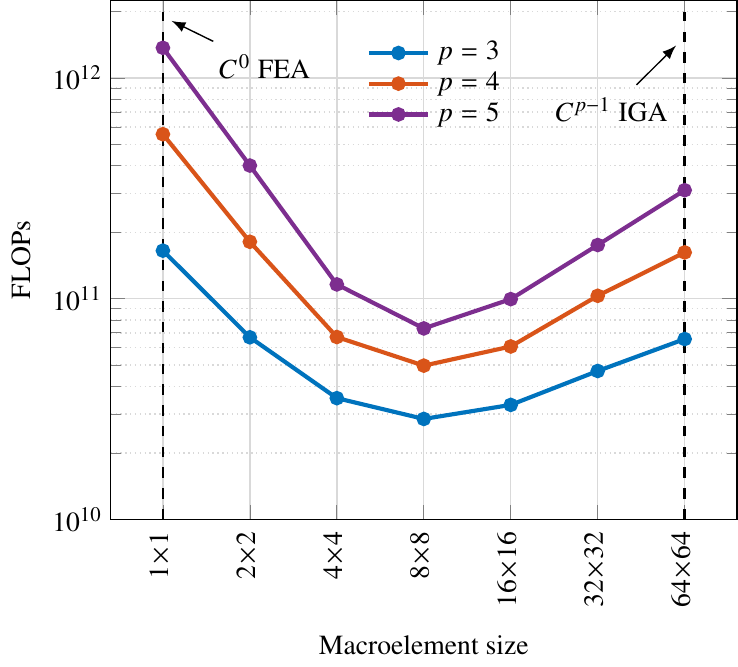}
    \caption{FLOPs ($64\times64$ elements)}
  \end{subfigure}\hspace{7pt}
  \begin{subfigure}{0.49\textwidth}\centering
    \pgfplotsset{height=0.85\textwidth}
    \includegraphics{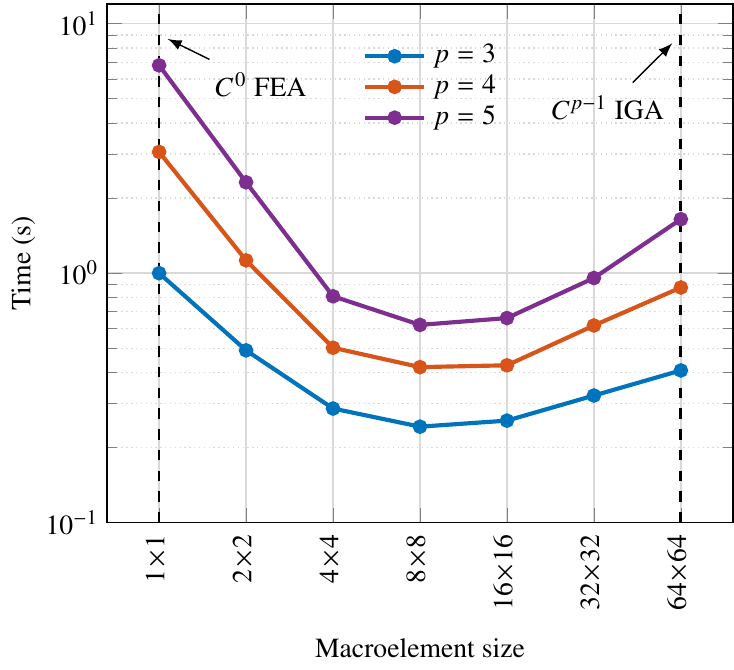}
    \caption{Time ($64\times64$ elements)}
  \end{subfigure} \\[6pt]
  \begin{subfigure}{0.49\textwidth}\centering
    \pgfplotsset{height=0.85\textwidth}
    \includegraphics{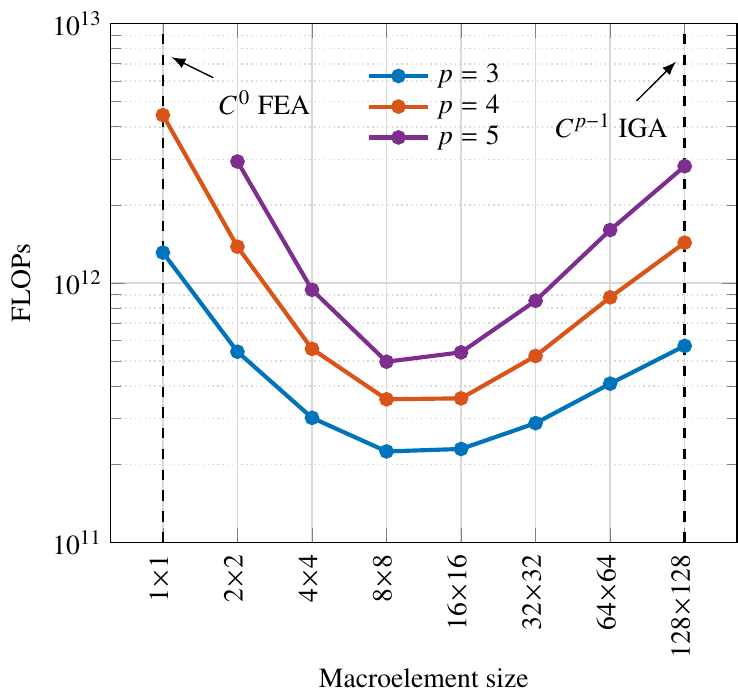}
    \caption{FLOPs ($128\times128$ elements)}
  \end{subfigure}\hspace{7pt}
  \begin{subfigure}{0.49\textwidth}\centering
    \pgfplotsset{height=0.85\textwidth}
    \includegraphics{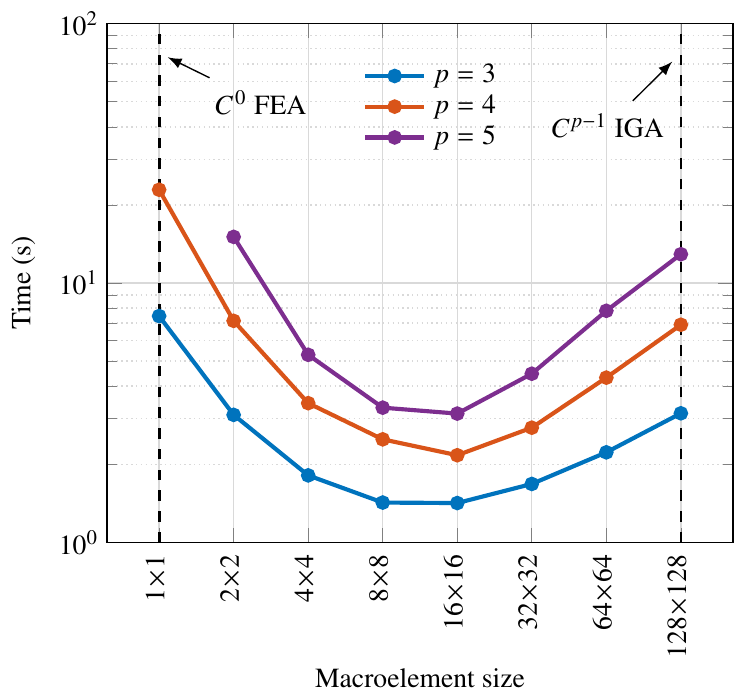}
    \caption{Time ($128\times128$ elements)}
  \end{subfigure}    
  \caption{Computational cost in terms of FLOPs and time for solving a 2.5D borehole resistivity problem per logging position per Fourier mode. 
  We test rIGA discretizations with two different grids of ${64\times64}$ and ${128\times128}$ elements.
  The computational times correspond to the use of parallel solver PARDISO using two threads.}
  \label{fg:TimeFLOPs}
\end{figure}

To numerically assess the computational efficiency
confirming the aforementioned theoretical results,
we consider two different grids in $\OM_{x,z}$ with ${64\times64}$ and ${128\times128}$ elements, respectively.
Using continuity reduction,
we split the mesh symmetrically into macroelements whose sizes are powers of two.
In this context, 
the maximum-continuity $C^{p-1}$ IGA discretization is composed of one macroelement containing the entire grid,
while FEA with $C^{\.0}$ continuity across all element interfaces is composed of macroelements that contain only one element.
\fig{\ref{fg:TimeFLOPs}} shows the number of FLOPs and time required to solve the borehole resistivity problem for each Fourier mode per logging position.
We compare the computational costs for different polynomial degrees and different continuity reduction levels of basis functions.
The cost of rIGA reaches the minimum with ${8\times8}$ macroelements almost in all cases, confirming the theoretical estimates obtained from the results of~\cite{Garcia2017}.

Our numerical tests show that for a moderate size 2.5D problem, the reduction in the number of FLOPs is $\Or(p)$ with respect to IGA. 
When compared to FEA, rIGA delivers larger improvement factors.
These improvement factors in terms of FLOPs also hold in terms of time when performing a sequential factorization.
In our parallel PARDISO solver, we observe a small degradation of the rIGA improvement factors in terms of times in comparison to those obtained in terms of FLOPs
(see \tab{\ref{tb:TimeFLOPs}}).

\begin{table}[!h]
\centering
\caption{Computational cost for the 2.5D borehole resistivity measurements per logging position per Fourier mode. We report the solution time and FLOPs when using $C^{p-1}$ IGA, rIGA with ${8\times8}$ macroelements, and $C^{\.0}$ FEA with the same number of elements and polynomial degree. The computational times correspond to the use of parallel solver PARDISO using two threads.}
\label{tb:TimeFLOPs}
\small
\begin{tabular}{@{}lllllllll@{}}
\toprule
\begin{tabular}[c]{@{}l@{}}Domain\\size\end{tabular} & 
\begin{tabular}[c]{@{}l@{}}Polynomial\\degree\end{tabular} & 
\begin{tabular}[c]{@{}l@{}}Discretization\\method\end{tabular} & 
\begin{tabular}[c]{@{}l@{}}Number of\\FLOPs\end{tabular} & 
\multicolumn{2}{l}{\begin{tabular}[c]{@{}l@{}}Improvement factor\\(FLOPs)\end{tabular}} & 
\begin{tabular}[c]{@{}l@{}}Time\\(s)\end{tabular}  & 
\multicolumn{2}{l}{\begin{tabular}[c]{@{}l@{}}Improvement factor\\(time)\end{tabular}} \\ \midrule
\multirow{9}{*}{64\.$\times$\.64} 
 & \multirow{3}{*}{3} 
 &	  IGA 	& 6.56e+10 & 		  &		 & 0.407 & 		  	& 		\\
 &  & rIGA 	& 2.85e+10 & IGA/rIGA & 2.30 & 0.242 & IGA/rIGA & 1.68 	\\
 &  & FEA 	& 1.65e+11 & FEA/rIGA & 5.79 & 0.999 & FEA/rIGA & 4.12 	\\[5pt]
 & \multirow{3}{*}{4} 
 & 	  IGA  & 1.62e+11 & 		 & 		 & 0.875 & 	   	  	& 		\\
 &  & rIGA & 4.97e+10 & IGA/rIGA & 3.26  & 0.419 & IGA/rIGA & 2.09  \\
 &  & FEA  & 5.56e+11 & FEA/rIGA & 11.19 & 3.061 & FEA/rIGA & 7.29  \\[5pt]
 & \multirow{3}{*}{5} 
 & 	  IGA & 3.10e+11 & 			& 		 & 1.645 & 		 	& 		\\
 &  & rGA & 7.33e+10 & IGA/rIGA & 4.23   & 0.620 & IGA/rIGA & 2.65  \\
 &  & FEA & 1.37e+12 & FEA/rIGA & 18.69  & 6.806 & FEA/rIGA & 10.98 \\ \cmidrule(){1-9}
\multirow{9}{*}{128\.$\times$\.128} 
 & \multirow{3}{*}{3} 
 & 	  IGA  & 5.72e+11 & 		 & 		& 3.144 & 		   & 		\\
 &  & rIGA & 2.24e+11 & IGA/rIGA & 2.55 & 1.423 & IGA/rIGA & 2.21 	\\
 &  & FEA  & 1.31e+12 & FEA/rIGA & 5.85 & 7.456 & FEA/rIGA & 5.24 	\\[5pt]
 & \multirow{3}{*}{4} 
 & 	  IGA  & 1.43e+12 & 		 & 		  & 6.903  & 		  & 	 \\
 &  & rIGA & 3.56e+11 & IGA/rIGA & 4.02   & 2.495  & IGA/rIGA & 2.77 \\
 &  & FEA  & 4.44e+12 & FEA/rIGA & 12.47  & 22.885 & FEA/rIGA & 9.17 \\[5pt]
 & \multirow{3}{*}{5} 
 & 	  IGA & 2.82e+12 & 			& 		 & 12.911 & 		 & 			\\
 &  & rGA & 4.97e+11 & IGA/rIGA & 5.67 	 & 3.305  & IGA/rIGA & 3.91 	\\
 &  & FEA & -- 		 & FEA/rIGA & -- 	 & -- 	  & FEA/rIGA & -- 		\\ \bottomrule
\end{tabular}
\end{table}


\subsection{Heterogeneous media} 
\label{sub:HeterogeneousMedia}

We further examine the accuracy of our rIGA approximation over two synthetic heterogeneous model problems.

\subsubsection{One geological fault}
\label{ssub:model1}

\begin{figure}[!h]
	\centering
	\includegraphics{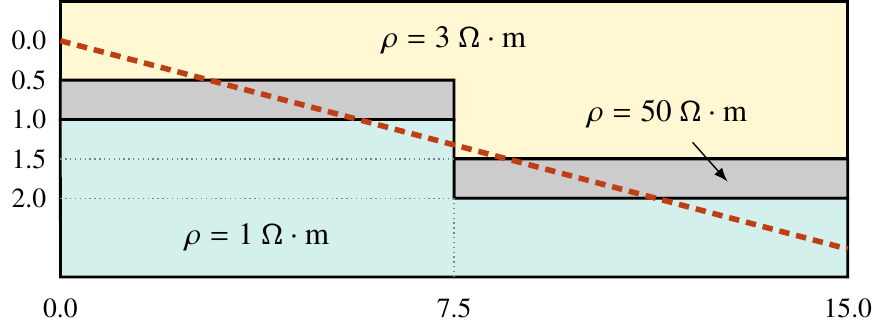}
	\caption{Model problem with a constant-angle well trajectory (red dashed line) passing through a geological fault and three different materials. Dimensions are in meters.}
	\label{fg:model1}
\end{figure}

We consider the model problem of \fig{\ref{fg:model1}}
with a constant dip angle of $80^\circ$.
We consider the LWD instrument described in \Sec{\ref{sec:ProblemDescription}}
and simulate measurements recorded over 200 equally-spaced logging positions throughout the well trajectory.

\begin{figure}[!h]
  \centering
  \pgfplotsset{width=0.8\textwidth,height=0.35\textwidth}
  \includegraphics{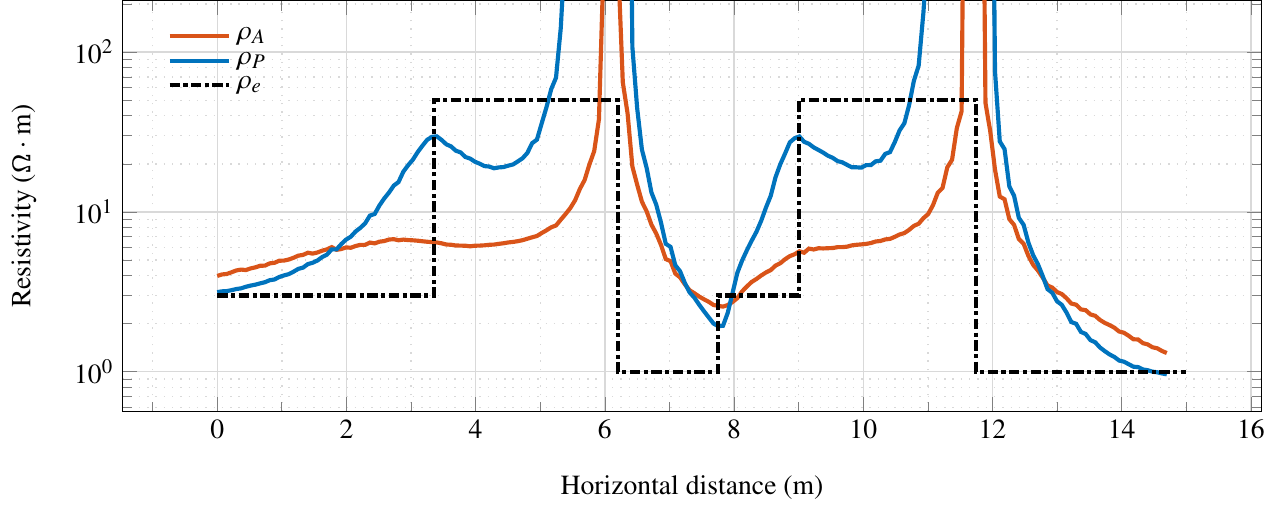}
  \caption{Apparent resistivities based on the attenuation, $\rho_A$, and phase, $\rho_P$, for the first model problem, compared with the real resistivity, $\rho_e$.
  We obtain the results using an rIGA discretization with ${64\times64}$ elements, ${p=4}$, and ${8\times8}$ macroelements.}
  \label{fg:AppRes1}
\end{figure}

\fig{\ref{fg:AppRes1}} shows the apparent resistivities based on the attenuation and phase (\,$\rho_A$ and $\rho_P$, respectively).
We obtain the results using ${N_f=70}$
and an rIGA discretization with ${64\times64}$ elements, ${p=4}$, and ${8\times8}$ macroelements. 
Results are in good agreement with those presented in~\cite{Rodrguez2018}.

\subsubsection{Two geological faults and inclined layers}
\label{ssub:model2}

\begin{figure}[!h]
	\centering
	\includegraphics{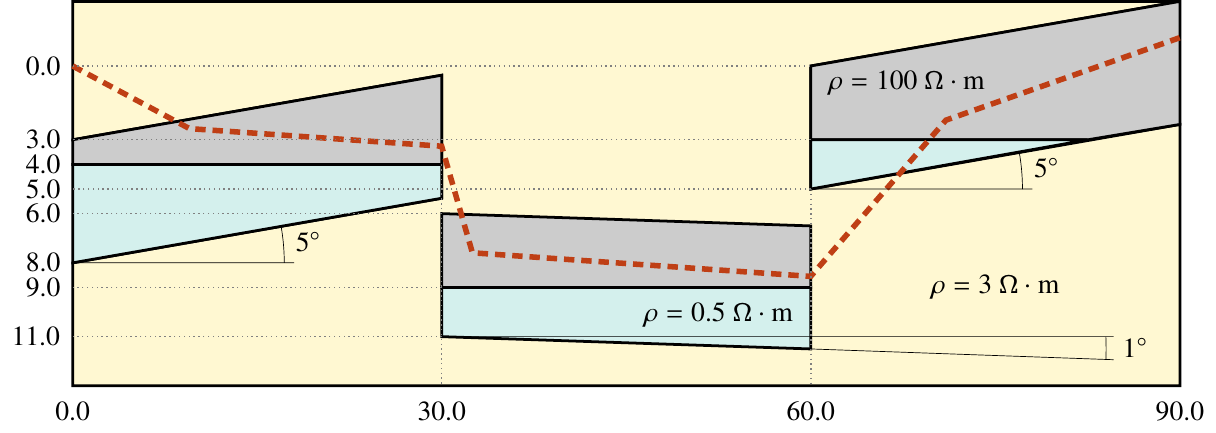}
	\caption{Second model problem with two geological faults and inclined layers. The tool trajectory (red dashed line) has different dip angles and passes through sandstone (yellow), oil-saturated (gray) and water-saturated (green) layers. Dimensions are in meters.}
	\label{fg:model2}
\end{figure}

\fig{\ref{fg:model2}} shows the second model problem containing two geological faults and inclined layers.
The logging trajectory starts from a sandstone layer with a resistivity of ${\rho=3~\Om\cdot\rm m}$, and passes through an oil-saturated layer with ${\rho=100~\Om\cdot\rm m}$. 
The tool trajectory also passes through a water-saturated layer with ${\rho=0.5~\Om\cdot\rm m}$.

In particular, inclined layers produce the so-called {\em staircase} approximations~\cite{Cangellaris1991}.
This phenomenon occurs because the physical interfaces of the conductivity model are not aligned with the element
edges.
Thus, the conductivity parameter takes different values inside some elements of the mesh. 
To tackle this issue, discretization techniques using nonfitting grids~\cite{Chaumont2018,Chaumont20182} are available, but they have not been considered here for simplicity.

\begin{figure}[!h]
	\centering
	\pgfplotsset{width=0.8\textwidth,height=0.4\textwidth}
	\includegraphics{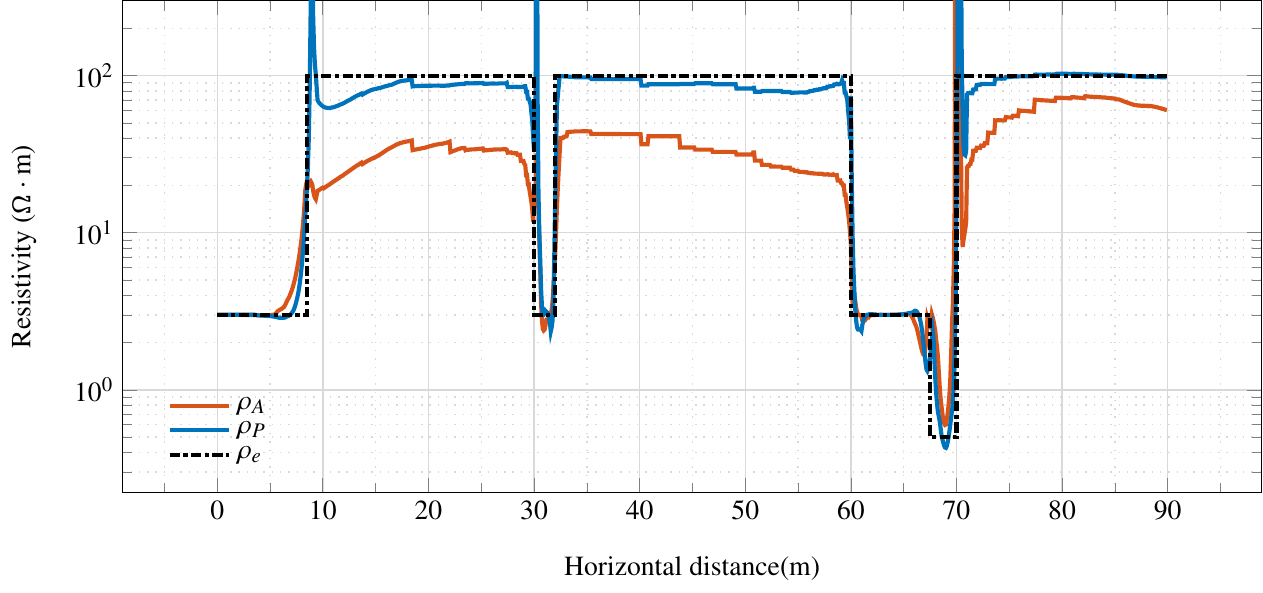}
	\caption{Apparent resistivities based on the attenuation, $\rho_A$, and phase, $\rho_P$, for the second model problem, compared with the real resistivity, $\rho_e$. We obtain the results using an rIGA discretization with ${64\times64}$ elements, ${p=4}$, and ${8\times8}$ macroelements.}
	\label{fg:AppRes2}
\end{figure}

\fig{\ref{fg:AppRes2}} shows the apparent resistivities based on the attenuation and phase throughout the logging trajectory and compares their value with the exact resistivity.
We simulate the resistivities at 1,080 logging positions with ${N_f=70}$. 
We use an rIGA discretization with ${64\times64}$ elements, ${p=4}$, and ${8\times8}$ macroelements.

\subsection{Database generation for DL inversion} 
\label{sub:Dataset}

\begin{figure}[!h]
	\centering
	\includegraphics{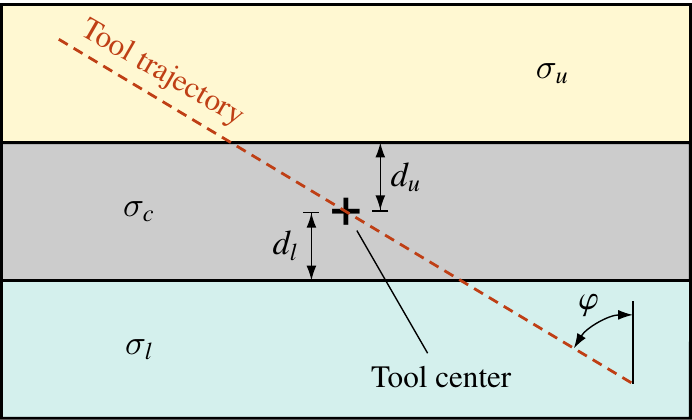}
	\caption{Varying parameters at each logging position when producing the training dataset for DL inversion.}
	\label{fg:DLmodel}
\end{figure}

To produce our synthetic training dataset for DL inversion, we consider heterogeneous medium containing three different layers and six varying parameters at each logging position, as described in \fig{\ref{fg:DLmodel}} and \tab{\ref{tb:DLmodel}}.
We select three different electrical conductivities: 
$\sg_c$ for the central layer, and $\sg_u$ and $\sg_l$ for the upper and lower layers, respectively.
We assume the tool center is always within the middle layer and has vertical distances of $d_u$ and $d_l$ from the upper and lower layers, respectively.
The sixth varying parameter is the dip angle, $\varphi$\., measured from the vertical direction.

In here, we create a dataset of 100,000 samples and compute the apparent resistivities obtained from
random combinations within a given range of resistivities ${\rho=1/\sg\in[1,100]~\Om\cdot\rm m}$ (see \tab{\ref{tb:DLmodel}}).

\begin{table}[!h]
\centering
\caption{Varying parameters employed to generate the training dataset for DL inversion.}
\label{tb:DLmodel}
\small
\begin{tabular}{@{}lll@{}}
\toprule
 Varying parameters                               &                          & \quad Interval              \\\midrule
 Electrical conductivity of the central layer     & \quad $\log_{10}(\sg_c$) & \quad $[-2,0]$          \\
 Electrical conductivity of the upper layer       & \quad $\log_{10}(\sg_u$) & \quad $[-2,0]$          \\
 Electrical conductivity of the lower layer       & \quad $\log_{10}(\sg_l$) & \quad $[-2,0]$          \\
 Distance of the tool center from the upper layer & \quad $\log_{10}(d_u$)   & \quad $[-2,1]$          \\
 Distance of the tool center from the lower layer & \quad $\log_{10}(d_l$)   & \quad $[-2,1]$          \\
 Dip angle between the tool and the layered media & \quad $\varphi$          & \quad $[80^\circ,100^\circ]$\\\bottomrule
\end{tabular}
\end{table}

For generating the dataset, we use two different types of parallelization.
One parallelization is related to the parallel factorization of the direct solver,
and the other is the trivial parallelization based on scheduling the solutions of independent earth models onto different processors. 
Using 40 threads,
we solve for 20 different earth models, each executing over two threads.
\tab{\ref{tb:TimeFLOPs}} shows that the required time for matrix factorization of the 2.5D EM problem using optimal rIGA discretization with ${64\times64}$ grid, ${p=4}$, and ${8\times8}$ macroelements is about 0.42 seconds per Fourier mode. 
Considering ${N_f=70}$, and the additional time required for pre/postprocessing and inter-thread communications, each set of independent runs (consists of 20 different earth models) takes about 40 seconds.
Thus, we perform 5,000 sequential runs to construct our 100,000 samples in about 56 hours. 
To create a larger database, we could execute over a cluster of hundreds of CPUs/threads, expecting a perfect parallel scalability.

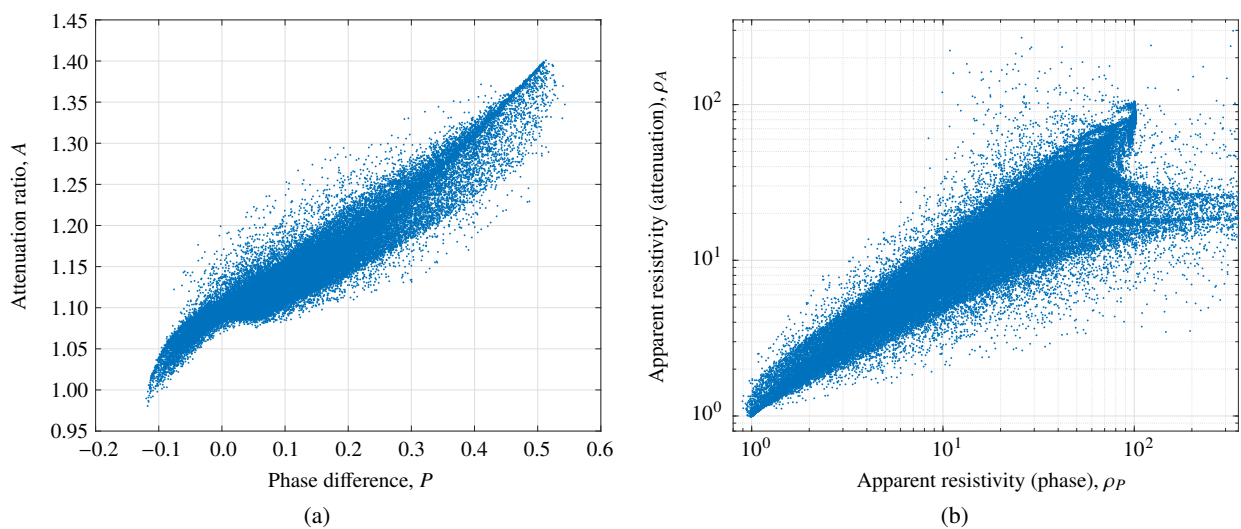
\begin{figure}[!h]
  \centering
  \begin{subfigure}{0.49\textwidth}\centering
    \begin{overpic}[width=1\linewidth,trim={0cm 0cm 2cm 0cm},clip]{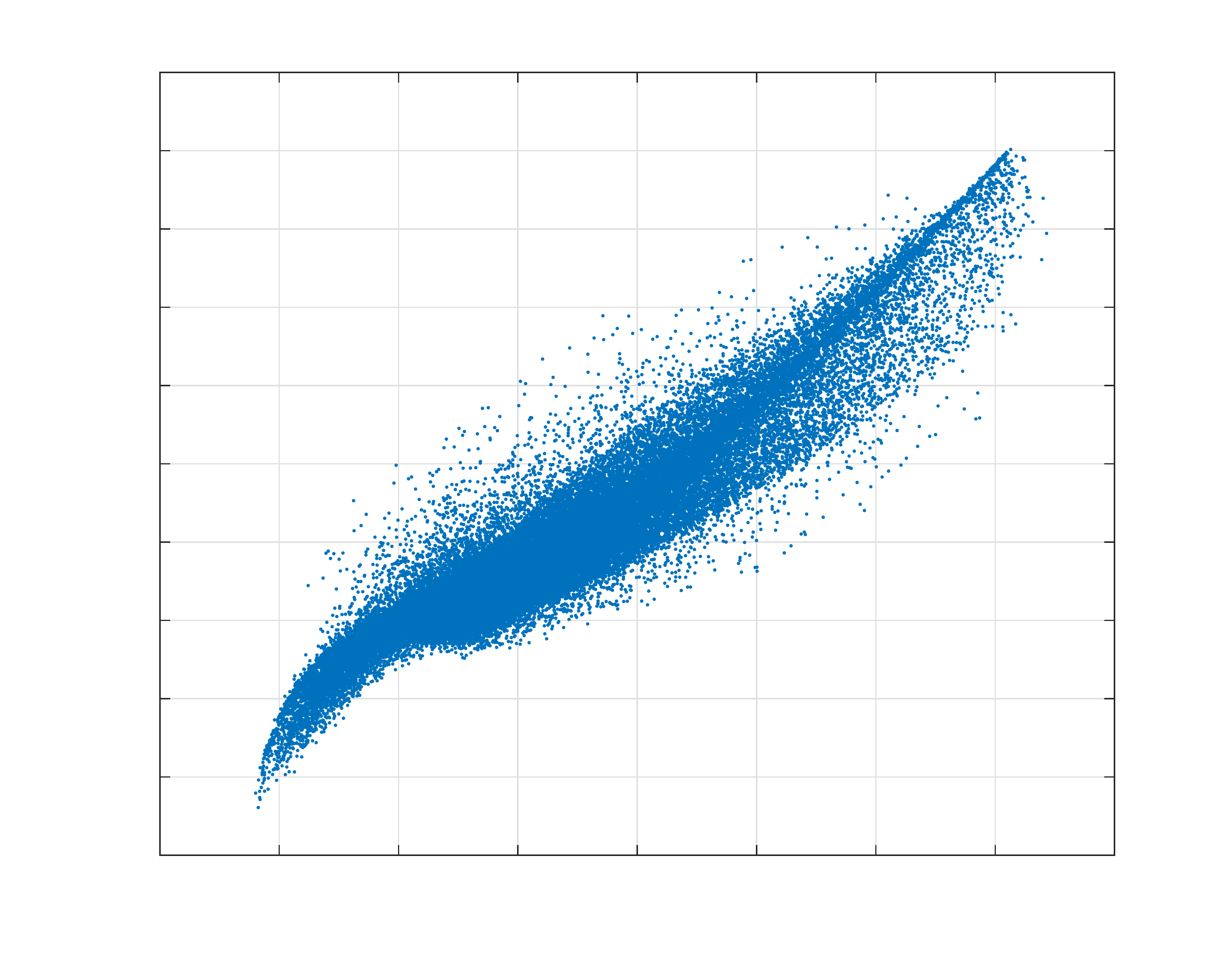}
      \put(42,0){\footnotesize Phase difference, $P$}
      \put(11.00,5){\footnotesize $-0.2$}
      \put(21.25,5){\footnotesize $-0.1$}
      \put(32.50,5){\footnotesize $ 0.0$}
      \put(42.75,5){\footnotesize $ 0.1$}
      \put(53.00,5){\footnotesize $ 0.2$}
      \put(63.25,5){\footnotesize $ 0.3$}
      \put(73.50,5){\footnotesize $ 0.4$}
      \put(83.75,5){\footnotesize $ 0.5$}
      \put(94.00,5){\footnotesize $ 0.6$}
      \put(6.5, 8.00){\footnotesize $0.95$}
      \put(6.5,14.75){\footnotesize $1.00$}
      \put(6.5,21.50){\footnotesize $1.05$}
      \put(6.5,28.25){\footnotesize $1.10$}
      \put(6.5,35.00){\footnotesize $1.15$}
      \put(6.5,41.75){\footnotesize $1.20$}
      \put(6.5,48.50){\footnotesize $1.25$}
      \put(6.5,55.25){\footnotesize $1.30$}
      \put(6.5,62.00){\footnotesize $1.35$}
      \put(6.5,68.75){\footnotesize $1.40$}
      \put(6.5,75.50){\footnotesize $1.45$}
      \begin{turn}{90} 
        \put(29,0){\footnotesize Attenuation ratio, $A$}
      \end{turn}
    \end{overpic}
    \caption{}
    \label{fg:DLmodelAppRes.a}
  \end{subfigure}\hspace{7pt}
  \begin{subfigure}{0.49\textwidth}\centering
    \begin{overpic}[width=1\linewidth,trim={0cm 0cm 2cm 0cm},clip]{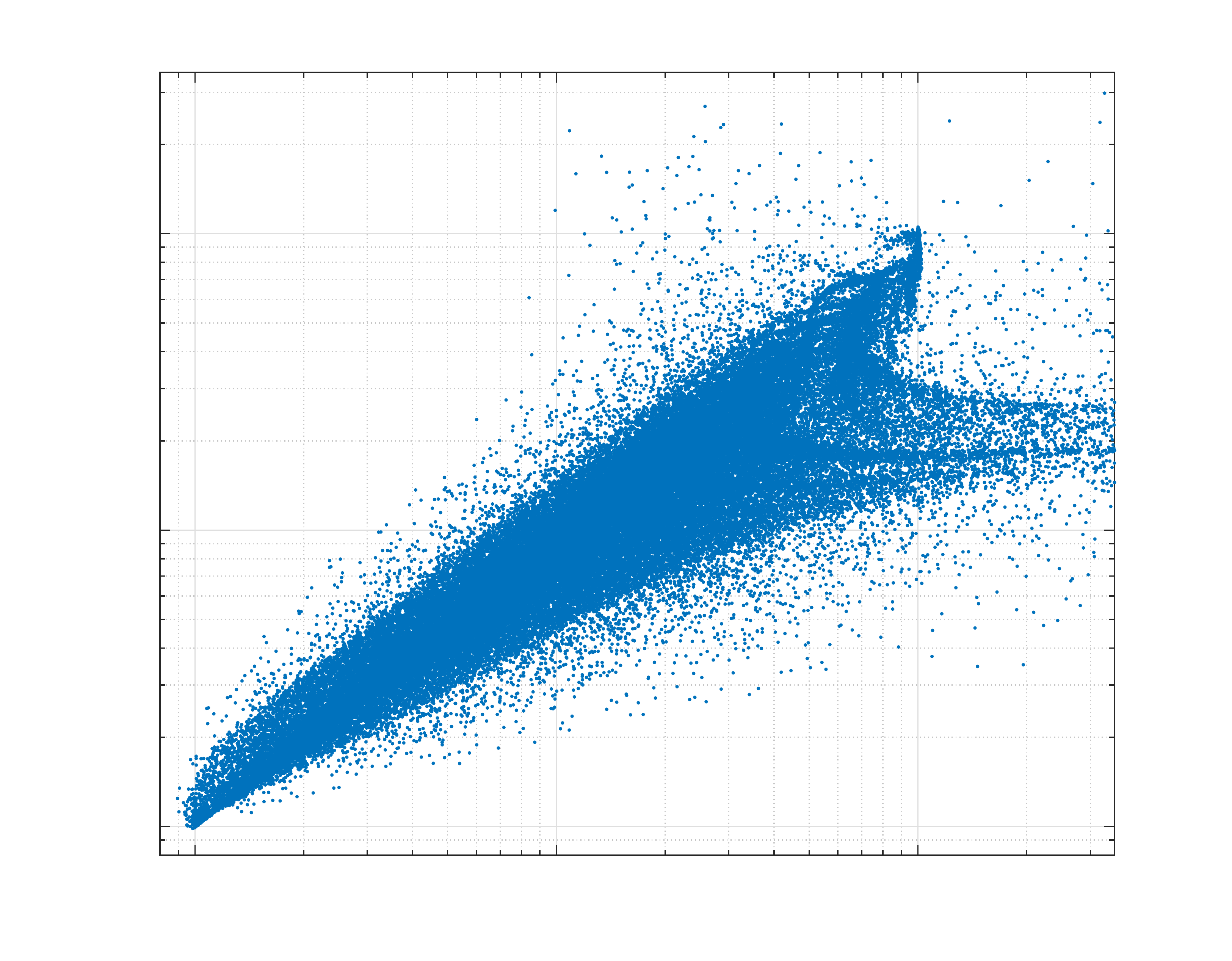}
      \put(34,0){\footnotesize Apparent resistivity (phase), $\rho_P$}
      \put(15.00,5){\footnotesize $10^{\.0}$}
      \put(46.25,5){\footnotesize $10^{\.1}$}
      \put(77.50,5){\footnotesize $10^{\.2}$}
      \put(7,11.00){\footnotesize $10^{\.0}$}
      \put(7,36.50){\footnotesize $10^{\.1}$}
      \put(7,62.00){\footnotesize $10^{\.2}$}
      \begin{turn}{90} 
        \put(17,0){\footnotesize Apparent resistivity (attenuation), $\rho_A$}
      \end{turn}
    \end{overpic}
    \caption{}
    \label{fg:DLmodelAppRes.b}
  \end{subfigure}
  \caption{(a) Attenuation ratio vs. phase difference, and (b) apparent resistivity based on attenuation vs. apparent resistivity based on phase,
  obtained for the 100,000 earth models. 
  We use rIGA discretization with ${64\times64}$ elements, ${p=4}$, and ${8\times8}$ macroelements for generating the database.}
  \label{fg:DLmodelAppRes}
\end{figure}

\fig{\ref{fg:DLmodelAppRes.a}} depicts the graphs of attenuation ratio, $A$, versus phase difference, $P$, obtained from the 100,000 earth models
when using rIGA discretization for generating the database.
Since there is a strong correlation between $A$ and $P$, the data distribution on the plot follows an almost straight line.
We also display in \fig{\ref{fg:DLmodelAppRes.b}} the correlation between apparent resistivities based on attenuation and phase.

\section{Conclusions} 
\label{sec.Conclusions}

We propose the use of refined isogeometric analysis (rIGA) discretizations for generating a (large) synthetic database for deep learning inversion of 2.5D borehole electromagnetic measurements. 
Such a database is essential for layer-by-layer estimation of the inverted earth models, which may be used for real-time adjustments of the well trajectory during geosteering operations. 

rIGA delivers computational savings of up to $\Or(p)$ compared to the high-continuity isogeometric analysis (IGA). 
When compared to a traditional finite element analysis (FEA) with the same mesh size and polynomial degree, rIGA provides higher improvement factors. 
At the same time, rIGA provides sufficiently accurate solutions for geosteering purposes.

To create a dataset for deep learning inversion,
we first selected certain discretization parameters based on the results of several homogeneous solutions.
Then, we checked the accuracy over homogeneous and heterogeneous media.
Finally,
we generated a meaningful synthetic database composed of 100,000 earth models with the corresponding measurements in about 56 hours using a workstation equipped with two CPUs.

As future work,
we propose the use of artificial intelligence based techniques to perform the inversion of a large set of borehole resistivity measurements with earth models containing multiple geological faults.
To create a larger database, we propose to execute the presented method over a cluster of hundreds of CPUs/threads.
We will also employ more complex earth model parameterizations including anisotropic layers.


\section*{Acknowledgment}

This work has received funding from 
the European Union's Horizon 2020 research and innovation program under the Marie Sklodowska-Curie grant agreement No 777778 (MATHROCKS), 
the European POCTEFA 2014-2020 Project PIXIL (EFA362/19) by the European Regional Development Fund (ERDF) through the Interreg V-A Spain-France-Andorra program, 
the Project of the Spanish Ministry of Science and Innovation with reference PID2019-108111RB-I00 (FEDER/AEI), 
the BCAM ``Severo Ochoa" accreditation of excellence (SEV-2017-0718), and the Basque Government through the BERC 2018-2021 program, 
the two Elkartek projects ArgIA (KK-2019-00068) and MATHEO (KK-2019-00085), 
the grant ``Artificial Intelligence in BCAM number EXP. 2019/00432", 
and the Consolidated Research Group MATHMODE (IT1294-19) given by the Department of Education.


\section*{References}
\bibliography{2.5D_refs}

\end{document}